\documentclass[conference]{IEEEtran}

  	\usepackage[pdftex]{graphicx}
  	\graphicspath{{../pdf/}{../jpeg/}}
	\DeclareGraphicsExtensions{.pdf,.jpeg,.png}
	\usepackage{algorithm}
	\usepackage{algpseudocode}
	\usepackage[cmex10]{amsmath}
	\usepackage{array}
	\usepackage{mdwmath}
	\usepackage{mdwtab}
	\usepackage{eqparbox}
	\usepackage{url}
\begin{document}
\title{\LARGE A case study of the Lunger phenomenon based on multiple algorithms}

 \author{\authorblockN{Zhang Jianan\authorrefmark{1}, Wang Yiyi\authorrefmark{2}, Duan Hongyi\authorrefmark{2}, Li Qingyang\authorrefmark{2} }
 \authorblockA{\authorrefmark{1}School of Mathematics, Shanghai University of Finance and Economics, Shanghai 200433, China}
 \authorblockA{\authorrefmark{2}Faculty of Electronic and Information, Xi'an Jiaotong University, Xi'an, 710049, China\\ziaqifei@ieee.org}}

\maketitle

\begin{abstract}
In this study, we conduct a thorough and meticulous examination of the Runge phenomenon. Initially, we engage in an extensive review of relevant literature, which aids in delineating the genesis and essence of the Runge phenomenon, along with an exploration of both conventional and contemporary algorithmic solutions. Subsequently, the paper delves into a diverse array of resolution methodologies, encompassing classical numerical approaches, regularization techniques, mock-Chebyshev interpolation, the TISI (Three-Interval Interpolation Strategy), external pseudo-constraint interpolation, and interpolation strategies predicated upon Singular Value Decomposition (SVD).

For each method, we not only introduce but also innovate a novel algorithm to effectively address the phenomenon. This paper executes detailed numerical computations for each method, employing visualization techniques to vividly illustrate the efficacy of various strategies in mitigating the Runge phenomenon. Our findings reveal that although traditional methods exhibit commendable performance in certain instances, novel approaches such as mock-Chebyshev interpolation and regularization-centric methods demonstrate marked superiority in specific contexts.

Moreover, the paper provides a critical analysis of these methodologies, specifically highlighting the constraints and potential avenues for enhancement in SVD decomposition-based interpolation strategies. In conclusion, we propose future research trajectories and underscore the imperative of further exploration into interpolation strategies, with an emphasis on their practical application validation. This article serves not only as a comprehensive resource on the Runge phenomenon for researchers but also offers pragmatic guidance for resolving real-world interpolation challenges.
\end{abstract}

\IEEEoverridecommandlockouts

\IEEEpeerreviewmaketitle

\section{Introduction and Review}

Polynomial interpolation is a pivotal issue in numerical analysis and scientific computation but can give rise to a well-known challenge: the Runge phenomenon. This phenomenon is characterized by an increase, rather than a decrease, in the error of polynomial interpolation over a finite domain with equidistant grids as the order of the interpolating polynomial increases \cite{she2019research}. This paper begins by reviewing various methods proposed to address this issue, including regularization techniques, Chebyshev sampling, least-squares fitting, methods based on radial basis functions, and interpolation methods incorporating external constraints.

Furthermore, this work explores the latest trends in deep learning to overcome the Runge phenomenon and the use of optimization techniques such as genetic algorithms. Regularization methods are a common approach to mitigate the Runge phenomenon. These methods aim to minimize a cost function comprising the sum of the residual and smoothness terms, maintaining the smoothness of the interpolation polynomial while fitting the data, and are widely applied in numerical analysis \cite{boyd1992defeating} \cite{jung2011simple}. One significant aspect of regularization methods is their adaptability to the specifics of a problem to achieve optimal performance in various scenarios.

Another effective approach to counteract the Runge phenomenon is Chebyshev node sampling. The phenomenon does not occur when interpolation points are sampled according to Chebyshev nodes. Hence, one method involves mimicking the Chebyshev approach by resampling equidistant interpolation points to approximate the distribution of a Chebyshev grid \cite{de2015constrained}. This method is advantageous in maintaining the smoothness of the interpolation polynomial while offering better numerical stability and is widely used in signal processing, image processing, and scientific computing.

Least squares fitting is a technique to reduce the Runge effect. It involves fitting \(n\) samples with a polynomial of degree \(d\), where \(d\) is significantly less than \(n-1\). This method is advantageous for better adapting to the noise and uncertainties in data, thus reducing interpolation errors \cite{boyd2009divergence}. Additionally, the least squares fitting method is versatile in handling non-equidistant interpolation points, making it more flexible for practical applications.

Radial basis function-based methods have been highly successful in overcoming the Runge phenomenon. These methods use radial basis functions to construct interpolation polynomials, known for their excellent approximation capabilities, especially in higher-dimensional spaces. Their application is extensive in interpolation problems, including geosciences, financial modeling, and medical image processing \cite{boyd2010six}.

Interpolation methods based on external constraints combine additional information to reduce interpolation errors. A common application of these methods is in image processing, where interpolation often accompanies constraints like boundary conditions or image features. These methods allow for the integration of these constraints into the interpolation process, thereby enhancing accuracy and stability \cite{belanger2017external}.

Recently, researchers in the field of deep learning have begun exploring the use of optimization algorithms to overcome the Runge phenomenon. They focus on issues associated with equidistant sampling points using parametric curve interpolation, searching for globally optimal parameter sequences. The potential for these methods in deep learning is vast, not only for solving the Runge phenomenon but also for other interpolation and fitting problems. Genetic algorithms (GAs), based on biological evolutionary principles, are widely used to solve various complex non-linear optimization challenges. GAs evolve a set of solutions (a population) through operations such as selection, crossover, and mutation to find the optimal solution \cite{gallagher1994genetic}. One study introduced an IGA-based optimal parameter search algorithm to counteract the Runge phenomenon \cite{lin2015searching}. This method leverages the diverse search and optimization capabilities of genetic algorithms to find parameters suitable for mitigating the Runge phenomenon.

Overcoming the Runge phenomenon has been a significant issue in polynomial interpolation. This paper reviews and replicates various methods, including regularization, Chebyshev sampling, least squares fitting, radial basis function-based methods, and interpolation methods based on external constraints. Each of these methods excels in different application domains, providing potent tools for addressing the Runge phenomenon. Future research will explore further improvements and applications of these methods to tackle a broader range of interpolation and fitting problems. The continued study of overcoming the Runge phenomenon is expected to propel advancements in numerical analysis and scientific computing.

\section{Visualization of Runge Phenomenon}
The Runge phenomenon is characterized by oscillations in the interpolating polynomial near the interval endpoints as the polynomial's order increases. This issue stems from the instability of higher-order polynomials at these endpoints. To illustrate this phenomenon, consider the function defined as:
\[
f(x) = \frac{1}{1 + 25x^2}
\]

\begin{figure}[ht!] 
\centering
\includegraphics[width=3.0in]{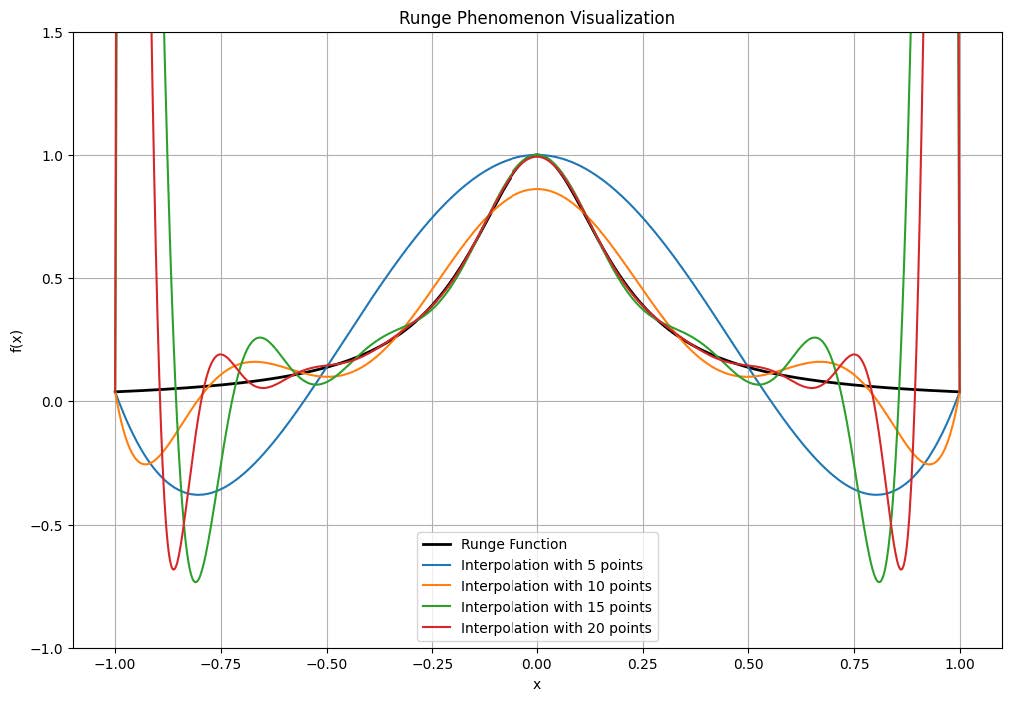}
\caption{Visualization of the Runge Phenomenon}
\label{1}
\end{figure}

Figure \ref{1} depicted by the black line, illustrates the original Runge function, which maintains smoothness throughout the interval. A noticeable trend emerges as the number of interpolation points is incremented: the polynomial interpolation's oscillations near the endpoints become increasingly pronounced. Although interpolation in the central region appears more aligned with the true function, the endpoint oscillations intensify. Specifically, with 5 interpolation points, the polynomial approximation closely mirrors the original function, yet slight endpoint oscillations are discernible. With 10 points, these oscillations near the endpoints are more pronounced. At 15 and 20 points, despite the central interval interpolation appearing more precise, the endpoint oscillations are markedly severe.

\section{Algorithmic principles and numerical computation reproduction}
\subsection{Traditional numerical methods: Chebyshev node-based interpolation and cubic spline interpolation}
\subsubsection{Chebyshev Nodes Interpolation}

For the function \( f(x) \), we consider using the first-kind Chebyshev nodes for interpolation within the interval \([-1, 1]\) to mitigate the Runge phenomenon. The first-kind Chebyshev nodes are defined as follows:
\begin{equation}
T_n(x) = \cos(n \arccos(x)), \quad x \in [-1, 1]
\end{equation}
where \( n \) is the degree of the polynomial.

The interpolation nodes on the interval \([-1,1]\) are the roots of the Chebyshev polynomial of the first kind:
\begin{equation}
x_k = \cos\left(\frac{(2k + 1)\pi}{2n + 2}\right), \quad k = 0, 1, \ldots, n
\end{equation}
This selection of nodes in the interval can effectively suppress the increase of interpolation error at the interval ends, which is the characteristic of the Runge phenomenon.

During the Chebyshev interpolation process, the error between the actual function and the interpolation polynomial will tend to be uniformly distributed across the entire interval rather than accumulating at the ends. This can effectively suppress the error magnification typically observed in the Runge phenomenon. The error bound is given by:
\begin{equation}
\| f(x) - I_n(x) \|_{\infty} \leq \frac{M_{n+1}}{(n+1)!} \| w_{n+1}(x) \|_{\infty} = \frac{M_{n+1}}{2^n(n + 1)}
\end{equation}

\subsubsection{Cubic Spline Interpolation}

Cubic spline interpolation utilizes piecewise cubic polynomials for segments delineated by interpolation nodes, ensuring match and twice differentiability at each node for enhanced function smoothness. The cubic splines \( S_i(x) \) within \( [x_i, x_{i+1}] \) are defined by:
\begin{itemize}
    \item \( S_i(x_i) = f(x_i) \), \( S_i(x_{i+1}) = f(x_{i+1}) \).
    \item \( S'_i(x_i) = S'_{i-1}(x_i) \), \( S''_i(x_i) = S''_{i-1}(x_i) \).
    \item \( S''_i(x) \) remains continuous, omitting \( x_{i+1} \) from the \( i \)-th interval.
\end{itemize}
This continuity of derivatives yields a composite function that is cubic and bi-differentiable. Cubic splines outperform polynomial interpolation by mitigating endpoint oscillations, thereby providing a precise representation of functions with sharp variations.

\subsubsection{Numerical Calculations and Results 1}

\begin{figure}[ht!] 
\centering
\includegraphics[width=3.0in]{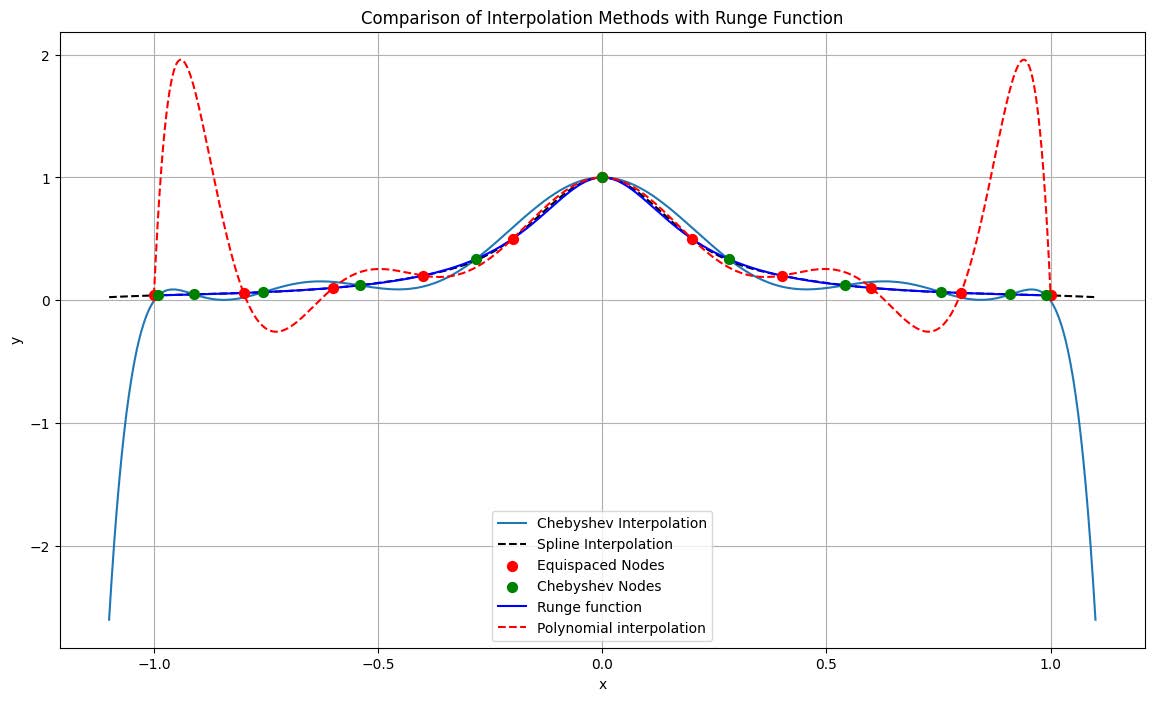}
\caption{Comparison of Interpolation Methods with Runge Function}
\label{2}
\end{figure}

Figure \ref{2} presents the experimental findings: The blue curve depicts the actual Runge Function, characterized by its centrosymmetry and peak near the interval's center. Chebyshev Interpolation closely approximates the function, underscoring its proficiency in addressing the Runge phenomenon with node clustering towards the interval ends to minimize interpolation error. Similarly, Cubic Spline Interpolation adheres closely to the central part of the function, with negligible deviations at the ends.

Node visualization: Uniformly spaced nodes are marked by red dots, while Chebyshev nodes are represented by green dots, with the latter's density increase near interval ends enhancing interpolation accuracy.

In essence, both Chebyshev and cubic spline interpolations effectively approximate the Runge function, particularly at the interval ends. They circumvent the Runge phenomenon inherent in high-order polynomial interpolations with uniformly spaced nodes. The results validate the exemplary performance of both methods, establishing their reliability for practical computational applications.

\subsection{Introduction of statistically based regularization methods}
\subsubsection{L1 Regularization methods}
In the context of L1 regularization, the Lasso model introduces an additional term, which imposes a constraint on the complexity of the model to prevent overfitting. The optimization problem can thus be articulated as:

\begin{equation}
    minimize \frac{1}{2N} \sum_{i=1}^{N} (y_i - w x_i)^2 + \alpha \sum_{j=1}^{p} |w_j|
\end{equation}

Where:

\begin{itemize}
    \item \( y_i \) represents the target variable,
    \item \( x_i \) denotes the feature variables,
    \item \( w \) is the model's weight vector,
    \item \( N \) is the number of observations,
    \item \( p \) indicates the number of features,
    \item \( \alpha \) is the regularization parameter, balancing the trade-off between fit and complexity.
\end{itemize}

The L1 norm of the weight vector encourages sparsity in the model coefficients, thereby promoting a model that integrates only the most significant predictors and enhancing the interpretability of the model. This is particularly beneficial when dealing with high-dimensional datasets where irrelevant features may dilute the predictive power of the model.

\subsubsection{L2 Regularization Methods}

L2 regularization, commonly associated with Ridge regression, contrasts L1 regularization by penalizing the square of the weights, which typically results in a model where the weight magnitudes are reduced uniformly. The optimization problem is formalized as follows:

\begin{equation}
    min \frac{1}{2N} \sum_{i=1}^{N} (y_i - w x_i)^2 + \alpha \sum_{j=1}^{p} w_j^2
\end{equation}

Unlike L1 regularization, L2 regularization does not result in sparsity of the model coefficients, thus it does not inherently perform feature selection, but rather it uniformly shrinks the coefficients to reduce model complexity.

\subsubsection{Elastic Net Regularization Methods}

Elastic Net regularization combines the properties of both L1 and L2 regularization. The objective function to be minimized is:

\begin{equation}
    min \frac{1}{2N} \sum_{i=1}^{N} (y_i - w x_i)^2 + \alpha \rho \sum_{j=1}^{p} |w_j| + \frac{\alpha (1-\rho)}{2} \sum_{j=1}^{p} w_j^2
\end{equation}

The parameter \( \rho \) governs the mix between L1 and L2 penalty terms. When \( \rho = 1 \), the Elastic Net is equivalent to Lasso regression; when \( \rho = 0 \), it becomes Ridge regression.

\subsubsection{Tikhonov Regularization}

Tikhonov regularization, also known as Ridge regression, stabilizes the inverse problems by introducing additional constraints. The method is particularly useful when the solution is expected to be smooth or when there is a need to stabilize the numerical inversion of a poorly conditioned matrix. The regularization modifies the least squares problem as follows:

\begin{equation}
    y = Ac
\end{equation}

Here, \( A \) is the Vandermonde matrix; \( c \) is the coefficient vector to be estimated. The Tikhonov regularization modifies the objective function to:

\begin{equation}
    c = \arg\min_c(\|Ac - y\|^2 + \|\Lambda c\|^2)
\end{equation}

Here, \( \Lambda \) is the regularization matrix, which typically is chosen to promote smoothness in the solution and represents the Tikhonov matrix.

\subsubsection{Numerical Calculations and Results 2.1}
This part of the numerical computation is performed for L1 regularization, L2 regularization, elastic network regularization

We evaluate the impact of various regularization techniques (No regularization, Ridge, Lasso, Elastic Net) on the Runge phenomenon. These techniques are tested using a set of Runge samples.

\noindent \textbf{Procedure and Results:}

\begin{enumerate}
    \item Data sampling.
    Using the function \texttt{np.linspace}, we generate an array of 11 equidistant points as the sample \( x \). A corresponding array of Runge function values is calculated as the sample \( y \), which serves as the basis for the comparison of different regularization techniques on the Runge phenomenon.

    \item Model training.
    Next, we train the models. The training data are augmented by adding polynomial features to increase the model complexity. The \texttt{polynomial\_regression} model is extended, incorporating Ridge regularization to mitigate the influence of higher-order polynomial terms, limiting their coefficients to a magnitude of 10. The results are summarized as follows:

    \begin{itemize}
        \item No regularization: Overfitting is evident in the model trained without regularization, resulting in an exaggerated depiction of the Runge phenomenon.

        \item L2 regularization (Ridge regression): The Ridge regularization method attenuates the Runge phenomenon, smoothing the model's behavior.

        \item L1 regularization (Lasso regression): Lasso regularization induces sparsity, simplifying the model by eliminating non-significant polynomial terms, thereby mitigating the Runge phenomenon.

        \item Elastic Net regularization (Elastic Net regression): Elastic Net regularization combines the features of both L1 and L2 regularization, effectively balancing between model complexity and sparsity, and thus, diminishes the Runge phenomenon's oscillations.
    \end{itemize}
\end{enumerate}

\begin{figure}[ht]
\centering
\includegraphics[width=0.5\textwidth]{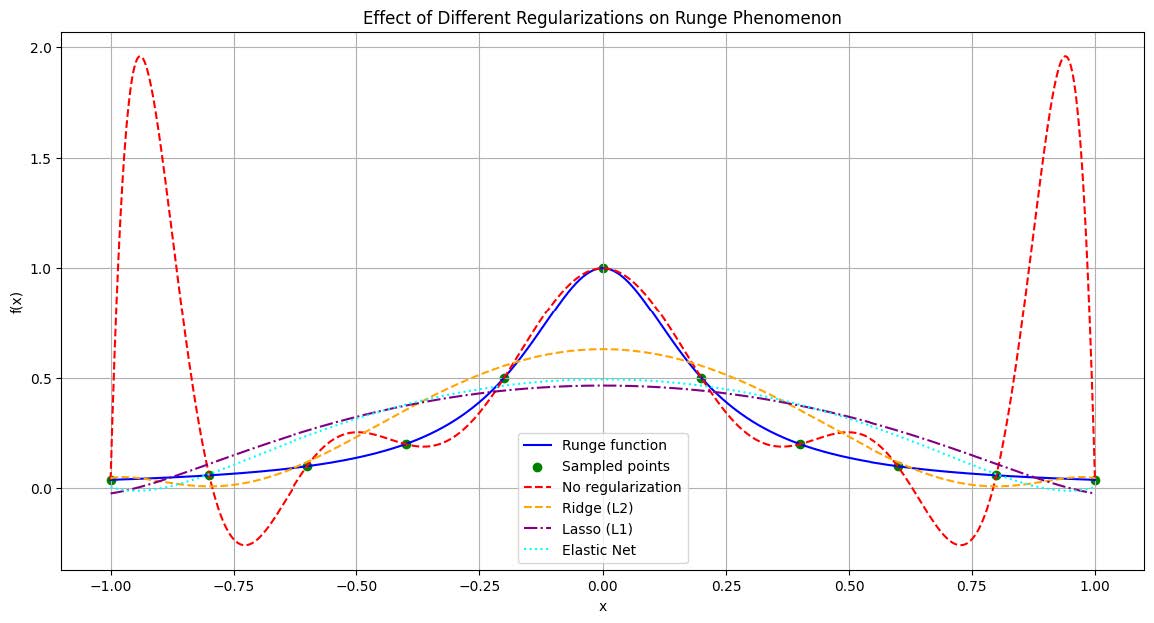}
\caption{Effect of Different Regularizations on Runge Phenomenon}
\label{3}
\end{figure}

The comparison outcomes are visually illustrated in Figure \ref{3}, providing a clear depiction of the effects of different regularization techniques on the Runge phenomenon.

\subsubsection{Numerical Calculations and Results 2.2}
This part of the numerical computation is performed for Tikhonov Regularization.

\textbf{Tikhonov Regularization:} Employ the \texttt{tikho\_inter} method, using polynomial degree expansion to model the data. Regularization is applied by using a lambda parameter \( \lambda \). The process can be summarized as follows:
\begin{itemize}
    \item Construct a Vandermonde matrix \( A \), representing each data point raised to the polynomial powers.
    \item Regularize the solution by applying Tikhonov regularization.
    \item Solve the least squares problem \( \text{Lstsqr} \) to obtain the model coefficients \( \text{coeffs} \).
\end{itemize}

\textbf{Data Sampling:} Create a uniform array of 11 sample points using \texttt{np.linspace} and calculate the corresponding Runge function values for these samples.

\textbf{Application of Tikhonov Regularization:} Implement Tikhonov regularization by fitting a high-degree polynomial to the data. Choose a polynomial degree (\( \text{degree} \)) greater than 10 and set the regularization parameter \( \lambda \) to 0.01.

\textbf{Generation of Dense Data Points:} Generate a dense grid of points using \texttt{np.dense} to interpolate and approximate the Runge function values with the true function \( y_{\text{true}} \) and the Tikhonov approximation \( y_{\text{approx}} \).

\textbf{Visualization:} Use Matplotlib to create a visualization that compares the original Runge function with the effects of Tikhonov regularization. This illustration should provide insights into the regularization's impact on the Runge phenomenon, including the suppression of oscillations due to the regularization parameter \( \lambda \).

\begin{figure}[ht]
\centering
\includegraphics[width=0.5\textwidth]{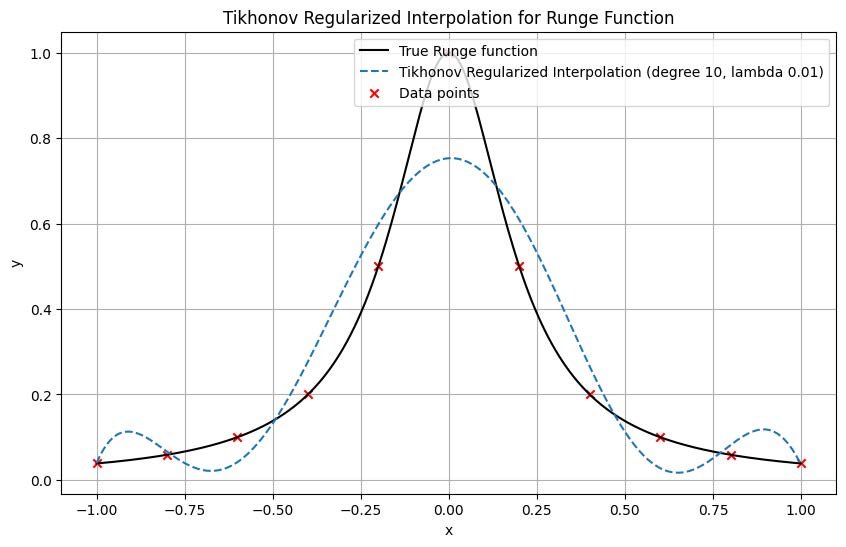}
\caption{Thikonov Regularize Interpolation for Runge Function}
\label{4}
\end{figure}

The results as shwon in Figure \ref{4} demonstrate that Tikhonov regularization effectively smoothens the Runge function's oscillatory behavior, proving to be a significant method for modeling phenomena with underlying regularization constraints.

\subsection{External Fake Constraints Interpolation}
\subsubsection{Overview of EFCI}

External Fake Constraints Interpolation (EFCI) introduces an innovative strategy to mitigate the Runge phenomenon by incorporating additional external constraints, often referred to as ``fake constraints,'' into the interpolation process. This method strategically utilizes constraints that are external to the original dataset.

The rationale behind EFCI is based on the Runge phenomenon, which occurs due to the equidistant spacing of interpolation points and leads to increased oscillations at the endpoints. By employing this novel technique, we can significantly reduce such undesired effects. EFCI imposes additional constraints on the interpolation polynomial \( P(x) \), ensuring that \( P''(x^*) = 0 \) where \( x^* \) are the external points, which results in the suppression of the polynomial’s curvature at its extremities. Consequently, the interpolation polynomial \( P(x) \) is coerced to align more closely with the function \( f(x) \) that it is intended to approximate, thus improving the overall interpolation accuracy.

The implementation of EFCI involves:

\begin{itemize}
    \item Constraint Specification: For each external point \( x_i \), we enforce:
    \begin{equation}
     P(x_i) = f(x_i) 
    \end{equation}

    \item Boundary Constraints: We introduce \( m/2 \) external points at each end of the interval \([-1, -1+\varepsilon]\) and \([1-\varepsilon, 1]\) to maintain the polynomial’s stability:
    \begin{equation}
    P''(x^*) = 0
    \end{equation}
\end{itemize}

Here, \( x^* \) denotes the external points.

By integrating these fake constraints into the interpolation process, the EFCI method effectively reduces the overshooting associated with the traditional interpolation over equidistant points without compromising the interpolation's fidelity.

\subsubsection{Numerical Calculations and Results 3}
In this section we will briefly describe and complete the numerical calculations based on EFCI.

\textbf{Definition of the Objective Function:} The objective function acts as the target for optimization problems, containing two error components. It starts by evaluating the polynomial values at specific x locations and then calculates the following error terms:
\begin{itemize}
    \item First error term: Function Value Matching. This measures the discrepancy between the polynomial values at the provided x locations and the true function values.
    \item Second error term: EFC Constraints. This assesses how the polynomial deviates at specific points (EFC positions) from the Runge function.
\end{itemize}

\textbf{Generation of Equidistant Nodes:} A series of equidistant x nodes is produced, and the Runge function values at these nodes are determined.

\textbf{Optimization: }The procedure iterates across different numbers of EFC (Equal Function Constraints) positions, from two to ten, creating EFC positions on either side of the function. An optimization challenge is then resolved to minimize the objective function, with the initial guesses for the polynomial coefficients set to zero and the x coordinates of the EFC positions.

\textbf{Optimal Results:} The experiment monitors and records the best optimization results, which are those with the smallest values of the objective function, along with the respective EFC positions and polynomial coefficients.

\textbf{Results Visualization:} A figure is created, displaying:
\begin{itemize}
    \item The genuine Runge function (solid blue line)
    \item EFC Interpolated Polynomial (dashed red line)
    \item Equidistant Nodes (green dots)
    \item EFC Positions (black crosses)
\end{itemize}

This figure vividly illustrates the performance of the interpolation polynomial with EFC against the original Runge function. The interpolation polynomial more accurately fits the shape of the Runge function around the EFC positions, diminishing oscillations.

\begin{figure}[ht]
\centering
\includegraphics[width=0.5\textwidth]{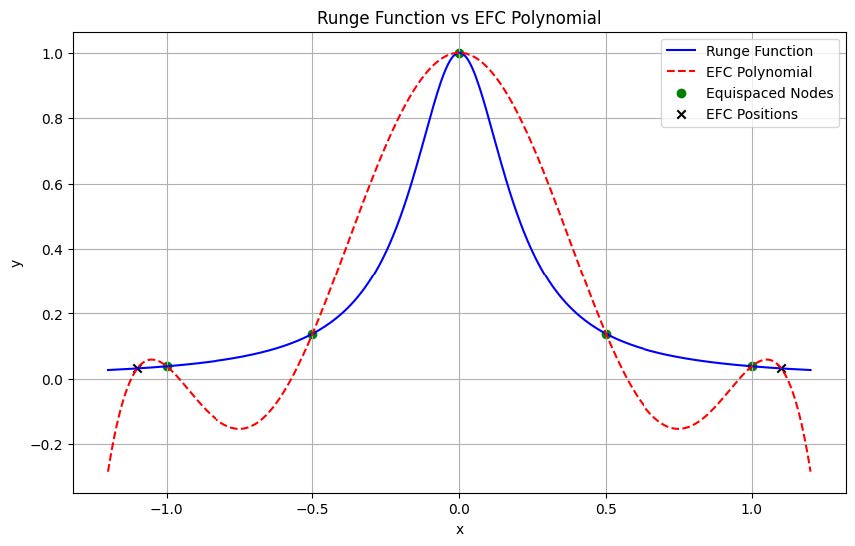}
\caption{Runge Function vs EFC Polynomial}
\label{5}
\end{figure}

The experimental results are shown in the Figure\ref{5}: During the experimental phase, by continuously varying the quantity of EFC positions, an optimal set of polynomial coefficients is found that reduces the oscillations in the Runge function interpolation. The graphical depiction provides a comparative view of the best-fit interpolation polynomial against the actual Runge function, emphasizing the influence of EFC positions. The optimal polynomial

\subsection{Mock-Chebyshev Method}
\subsubsection{General Principles}

The Mock-Chebyshev method is a strategy that approximates the optimal node distribution of the Chebyshev method, which is renowned for minimizing polynomial interpolation errors. This method selects nodes based on the extrema of Chebyshev polynomials, which tend to cluster more densely near the interval's endpoints, thereby mitigating the Runge phenomenon. However, obtaining exact Chebyshev nodes can be impractical in certain scenarios, especially when dealing with complex or unknown underlying functions. The Mock-Chebyshev method provides a surrogate by allocating pseudo-Chebyshev nodes that mimic the ideal distribution of actual Chebyshev nodes. These nodes, while not aligning with the true extrema of Chebyshev polynomials, offer a close approximation and effectively reduce interpolation errors.

Nodes are ascertained by approximating the Chebyshev node distribution within a given interval, although not achieving the exact Chebyshev polynomial extrema.

The Mock-Chebyshev method offers practical solutions when exact Chebyshev node distribution is unattainable, providing substantial benefits over the conventional equidistant node distributions. Notable advantages include:
\begin{itemize}
    \item Reduction of Error: It can diminish interpolation errors similar to the Chebyshev method, especially around the extremities of the interval.
    \item Flexibility: Given the challenges in acquiring precise Chebyshev nodes, this method allows for a flexible and adaptable approximation.
    \item Simplicity: It streamlines the node distribution process by eschewing the need for exact Chebyshev polynomial extrema calculations, adopting a mock distribution that closely emulates the Chebyshev ideal.
\end{itemize}

\subsubsection{Numerical Calculations and Results 4}

\textbf{Least-Squares Polynomial Fitting:} Initially, 20 uniformly distributed data points were generated within the interval \([-1, 1]\). The least-squares polynomial fitting method was then applied to these data points. This process was repeated 10 times to refine the fit.

\textbf{Mock-Chebyshev Subset Interpolation:} Afterward, a technique known as "Mock-Chebyshev Subset Interpolation" was employed. This method, based on the interaction between Chebyshev-Lobatto nodes and a full grid, involved selecting 11 Chebyshev-Lobatto nodes within the interval \([-1, 1]\), excluding the endpoints. A subset was created from the full grid's points closest to these Chebyshev-Lobatto nodes. The Mock-Chebyshev interpolation was then performed on this subset, iterating the polynomial fitting 10 times for optimal modeling of the data.

\textbf{Result Visualization:} Finally, the experimental results were visualized, showcasing the true Runge function curve (in blue), the least-squares fitting (dashed red line), and the Mock-Chebyshev interpolation results (dotted green line). The experimental data points (black circles) were also indicated.

A comparative analysis of the least-squares polynomial fitting and Mock-Chebyshev interpolation methods revealed their respective approximation accuracies to the Runge function. The differences between the methods and the actual Runge function are clearly displayed in the Figure\ref{6} below.

\begin{figure}[ht]
\centering
\includegraphics[width=0.5\textwidth]{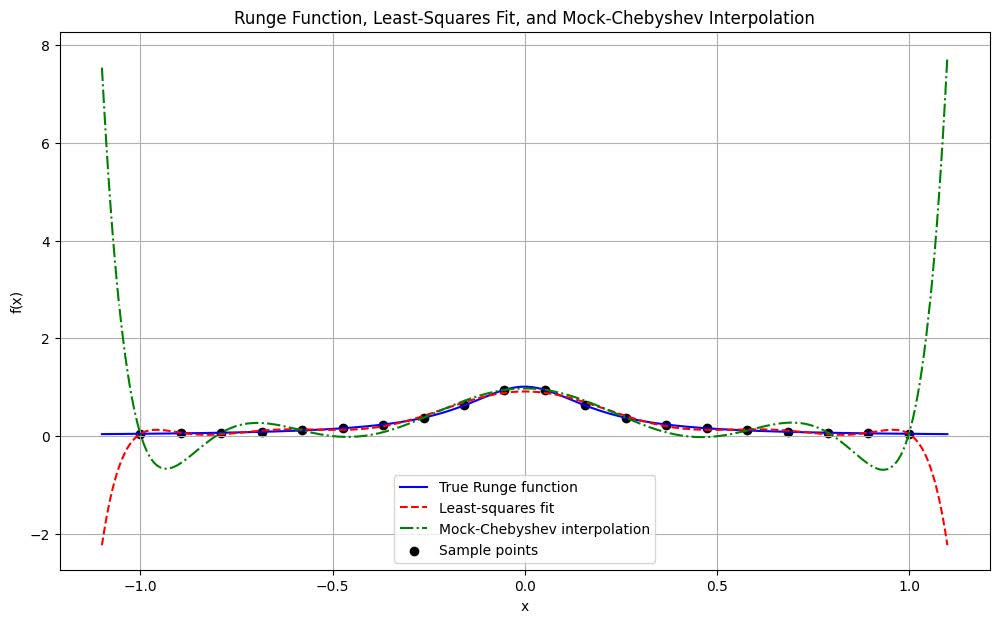}
\caption{Runge Function, Least-Squares Fit and Mock-Chebyshev Interpolation}
\label{6}
\end{figure}

\subsection{Three-Interval Smooth Interpolation (TISI)}
\subsubsection{TISI Algorithm Principles}

The TISI algorithm, standing for Three Interval Subdivision Interpolation, is an innovative method designed to enhance the interpolation of the Runge function. It subdivides the main interval into three subintervals, each with a tailored interpolation strategy.

\begin{itemize}
    \item \textbf{Left Interval:} In the left subinterval, interpolation proceeds in a traditional manner, as this area is less affected by the Runge phenomenon and therefore requires a standard interpolation approach.
    
    \item \textbf{Central Interval:} The central subinterval employs the Lagrange interpolation method, addressing the region's susceptibility to oscillations and ensuring stability and precision in the interpolation results.
    
    \item \textbf{Right Interval:} The right subinterval, significantly impacted by the Runge phenomenon, is approached with strategies such as:
    \begin{itemize}
        \item Increasing node density in regions with high oscillatory behavior.
        \item Modifying the interpolation method for x values within this interval, utilizing the Lagrange method.
        \item Applying strict Lagrange interpolation at the interval's boundaries to mitigate overshooting.
    \end{itemize}
\end{itemize}

These subinterval-specific strategies enable the TISI method to provide an accurate and stable interpolation across the entire interval, reflecting the TISI algorithm's principle of adapting interpolation techniques to subinterval characteristics for a smoother overall result.

\subsubsection{Numerical Calculations and Results 5}
For numerical simulations, the following procedures were adopted:

\begin{itemize}
    \item \textbf{Lagrange Interpolation for Function Fitting:} 
    
    The \texttt{lagrange\_interpolation} function performs the Lagrange interpolation, seeking a closer approximation to the true function within the influential interval.
    
    \item \textbf{Three-Interval Method Application:} 
    
    The \texttt{three\_interval\_method} conducts interpolation by dividing the domain into three parts, each adopting an appropriate method to minimize the impact of the Runge phenomenon.
    
    \item \textbf{Node Determination in Each Interval:} 
    
    In general consequences nodes are determined via the \texttt{three\_interval\_method} for each segment, tailoring the interpolation approach to the characteristics of each sub-domain.
    
    \item \textbf{Interval-Specific Node Selection:} 
   
    For each interval, a particular set of nodes (11 for each interval) is chosen, aiming to closely match the Runge function within that segment.
    
    \item \textbf{Selection of Appropriate Methods:} 
   
    Each interval may necessitate a different approach, and the \texttt{three\_interval\_method} facilitates the selection of the most suitable interpolation method for each interval.
\end{itemize}

\begin{figure}[ht]
\centering
\includegraphics[width=0.5\textwidth]{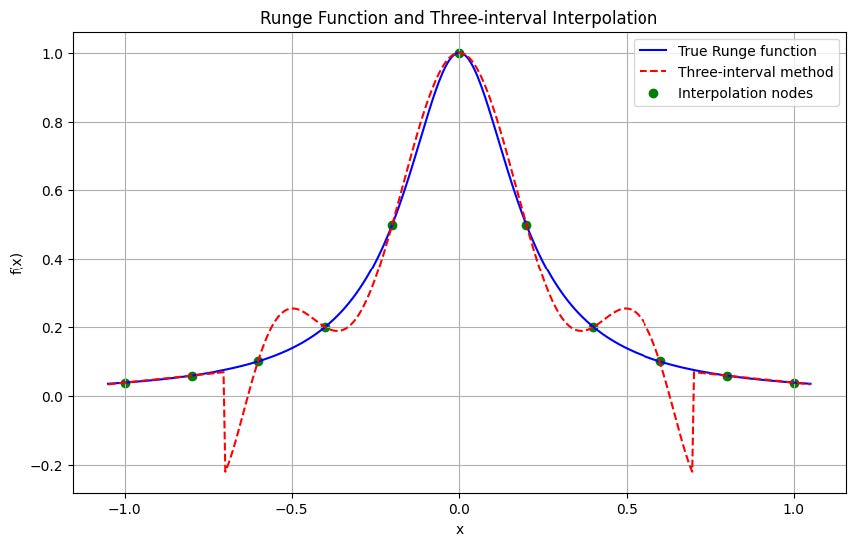}
\caption{Runge Function and Three-interval Interpolation}
\label{7}
\end{figure}

The results of the simulation are illustrated in Figure \ref{7}, which clearly depicts the efficacy of the three-interval method interpolation in reducing the oscillatory behavior associated with Runge's phenomenon.

\subsubsection{Improved Three-Interval Interpolation}

The improved three-interval method refines Runge function interpolation, particularly targeting the pronounced oscillations at the endpoints of the interval. It subdivides the domain into three specific subintervals, each with a customized approach:

\begin{itemize}
    \item \textbf{Left Interval:} Spanning from \( x = -1 \) to \( x = -1 + \epsilon \), where traditional interpolation is less challenged by the Runge phenomenon.
    
    \item \textbf{Right Interval:} Extending from \( x = 1 - \epsilon \) to \( x = 1 \), it adopts a symmetrical strategy to the left interval.

    \item \textbf{Central Interval:} Covering the middle portion, this area requires an enhanced interpolation method to address the significant oscillatory behavior due to the Runge phenomenon.
\end{itemize}

This methodological advancement, with its specialized interpolation strategies for each segment, secures a more stable and precise depiction of the Runge function. It is particularly adept at mitigating oscillations within the central subinterval, producing a smoother interpolation result.

The effectiveness of this refined approach is illustrated in Figure \ref{8} below.

\begin{figure}[ht]
\centering
\includegraphics[width=0.5\textwidth]{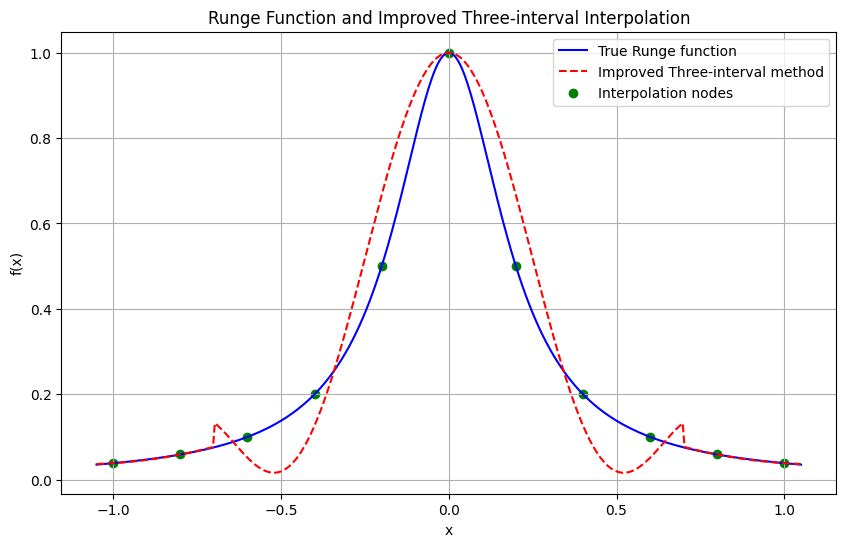}
\caption{Runge Function and Improved Three-interval Interpolation}
\label{8}
\end{figure}

\subsection{Restricted mock-Chebyshev least squares method}
\subsubsection{Methodology in Restricted mock-Chebyshev}

The methodology involves a series of advanced algorithms designed for the nuanced approximation of interpolation in functions with pronounced oscillations, such as the Runge function. The key steps include:

\begin{itemize}
    \item \textbf{Mock-Chebyshev Node Determination:} 
    The Mock-Chebyshev method determines a set of proxy nodes that closely approximate the optimal Chebyshev node distribution. These nodes provide a reliable base for improved interpolation, though they do not align exactly with the Chebyshev polynomial extrema.
    
    \item \textbf{Interval-Specific Node Selection:} 
    Nodes are selected within each interval with precision, tailored to the unique oscillatory characteristics of the function within those segments, thereby optimizing the interpolation accuracy.
    
    \item \textbf{Strategic Interpolation:} 
    The interpolation within the mock-Chebyshev framework is strategically adapted to the node distribution, particularly in regions with notable oscillations, to enhance stability and accuracy.
    
    \item \textbf{Iterative Refinement:} 
    An iterative process refines the node set within the mock-Chebyshev approach, aiming for an ideal node distribution that yields the most accurate interpolation of the target function.
\end{itemize}

Through these steps, the mock-Chebyshev method significantly improves the interpolation of functions that display strong oscillatory behaviors, effectively reducing the Runge phenomenon's typical errors.

\begin{figure}[ht]
\centering
\includegraphics[width=0.5\textwidth]{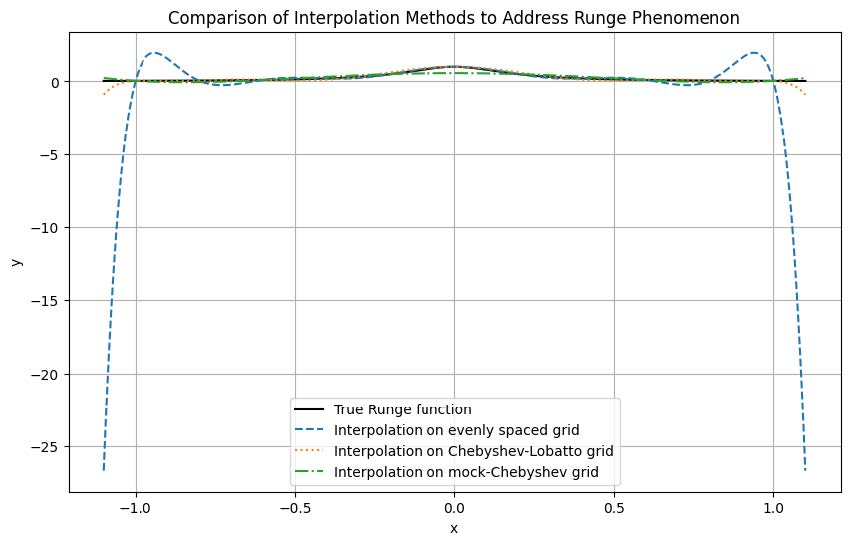}
\caption{Comparison of Interpolation Methods to Address Runge Phenomenon}
\label{9}
\end{figure}

\subsubsection{Numerical Calculations and Results 6}
In conjunction with other similar methods for comparison, the polynomial interpolation function is defined as follows: The \texttt{interpolate\_polynomial} function carries out polynomial interpolation based on a set of input x and corresponding y values.

Two different input data sets were generated, one on an evenly spaced grid and the other on a Chebyshev-Lobatto grid. These grids differ in their node selection methodology.

The Mock-Chebyshev subset, \texttt{x\_mock\_chebyshev} and \texttt{y\_mock\_chebyshev}, is derived by selecting every other point from the evenly spaced grid, emulating a Chebyshev subset.

Interpolation is performed on a dense \texttt{x\_dense} grid using various methods (evenly spaced, Chebyshev-Lobatto interval, Mock-Chebyshev subset), achieving corresponding function values.

The visualization of results includes:
\begin{itemize}
    \item The true Runge function (solid black line)
    \item Interpolation on an evenly spaced grid (dashed line)
    \item Interpolation on a Chebyshev-Lobatto interval (dotted line)
    \item Interpolation on a Mock-Chebyshev subset (dash-dot line)
\end{itemize}

As depicted in Figure \ref{9}, these representations distinctly highlight the differences in interpolation methods, especially in areas where the Runge function's oscillations are prominent, assessing the efficacy of each method.

\begin{figure}[ht]
\centering
\includegraphics[width=0.5\textwidth]{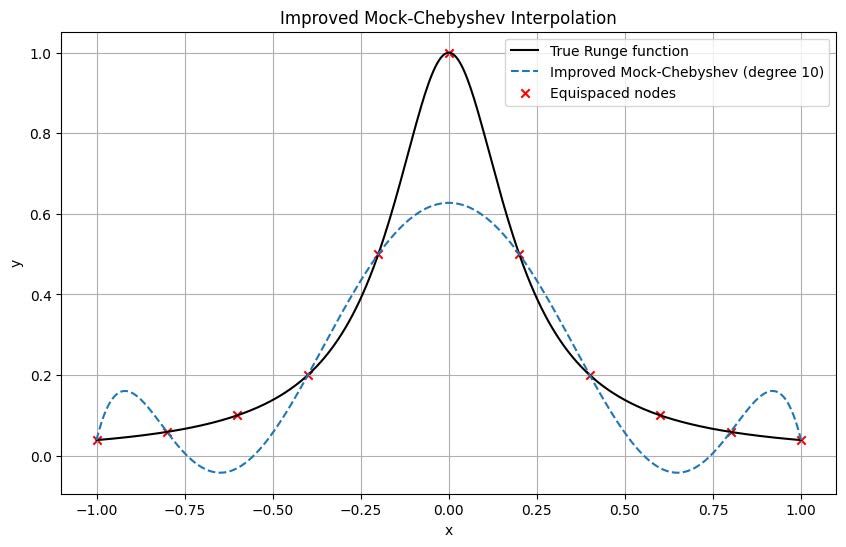}
\caption{Improved Mock-Chebyshev Interpolation}
\label{10}
\end{figure}

Figure \ref{10} illustrates the Improved Mock-Chebyshev Interpolation.

\subsection{SVD-Based Interpolation Method}

\subsubsection{Algorithm Principles}
The core idea of this approach utilizes the analytic properties of functions, applying Singular Value Decomposition (SVD) to the interpolation matrix and proceeding with appropriate truncation.
\begin{itemize}
    \item Begin by interpolating the Runge function.
    \item Construct the interpolation matrix using either Legendre polynomials or basic monomials.
    \item Apply SVD to the interpolation matrix.
    \item Truncate the singular values based on a predefined threshold.
    \item Reconstruct the interpolation matrix with the truncated singular values.
\end{itemize}

\subsubsection{Numerical Calculations and Results 7}
This section conducts an experiment with SVD truncation interpolation on a uniform grid to address the Runge phenomenon, evaluating the effect of different truncation thresholds on interpolation accuracy.
\begin{itemize}
    \item The \texttt{legendre\_matrix} function generates the interpolation matrix from Legendre polynomials, a fundamental tool for polynomial interpolation.
    \item The \texttt{svd\_truncated\_interpolation} function performs SVD truncation interpolation, taking an input data set 'x' with corresponding function values 'y', polynomial degree 'degree', and truncation threshold 'threshold'. It decomposes the interpolation matrix using SVD, truncates based on the threshold, and then reconstructs the polynomial coefficients for interpolation.
    \item Define a uniform grid with `x\_uniform` containing 11 nodes and `y\_uniform` holding the Runge function values at these nodes.
    \item Use different truncation thresholds listed in 'thresholds', such as 1e-2, 1e-5, 1e-10, and 1e-15, for comparative analysis of interpolation results.
    \item Generate a dense set of x values, `x\_fine`, for plotting detailed graphs.
\end{itemize}

\begin{figure}[ht]
    \centering
    \includegraphics[width=0.5\textwidth]{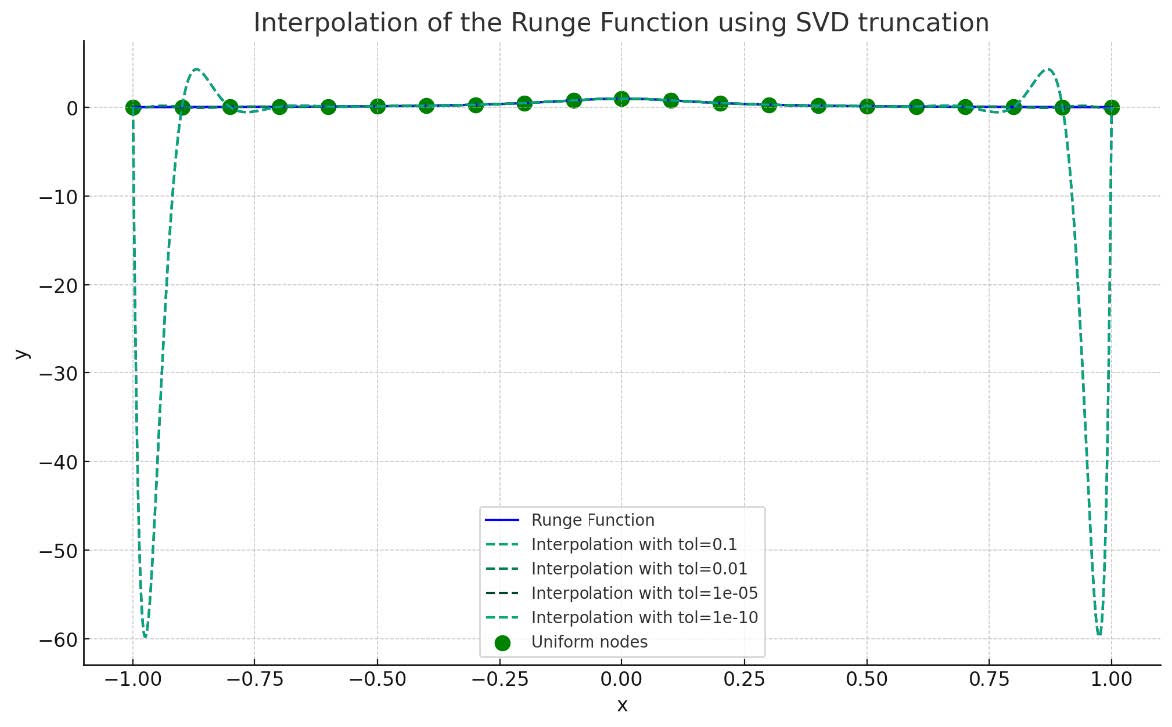}
    \caption{Interpolation of the Runge Function using SVD truncation.}
    \label{fig:svd_truncation}
\end{figure}

\begin{figure}[ht]
    \centering
    \includegraphics[width=0.5\textwidth]{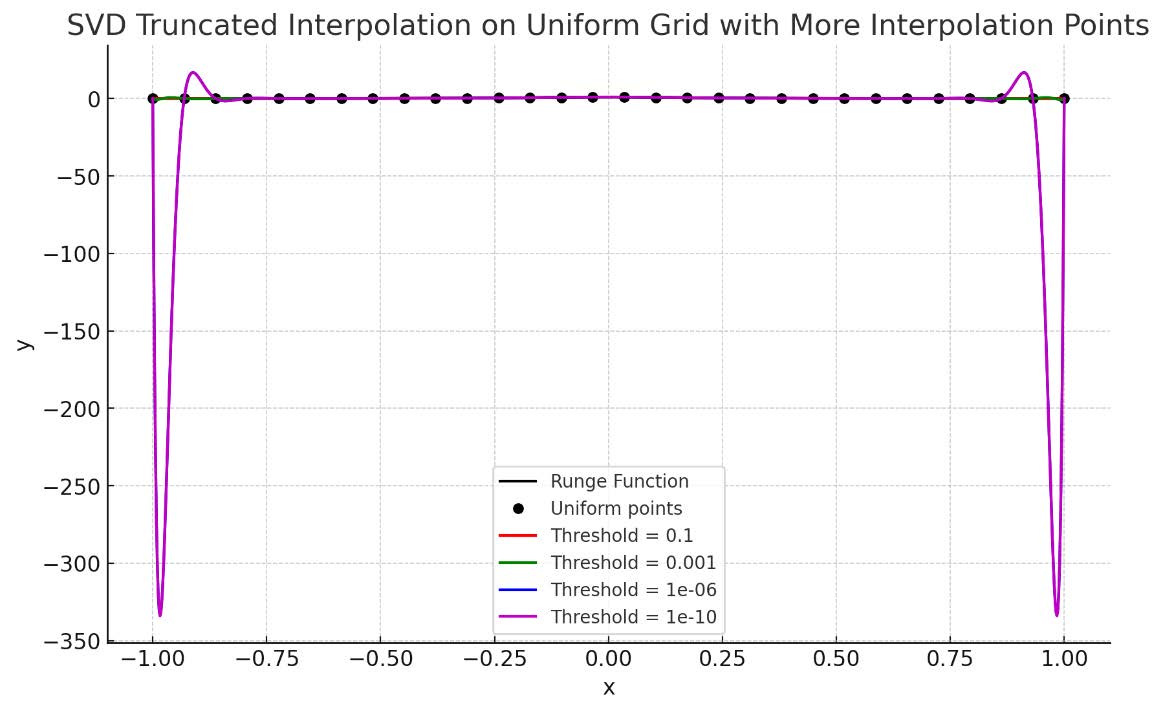}
    \caption{SVD Truncated Interpolation on Uniform Grid with More Interpolation Points.}
    \label{fig:svd_uniform_grid}
\end{figure}

\begin{figure}[ht]
    \centering
    \includegraphics[width=0.5\textwidth]{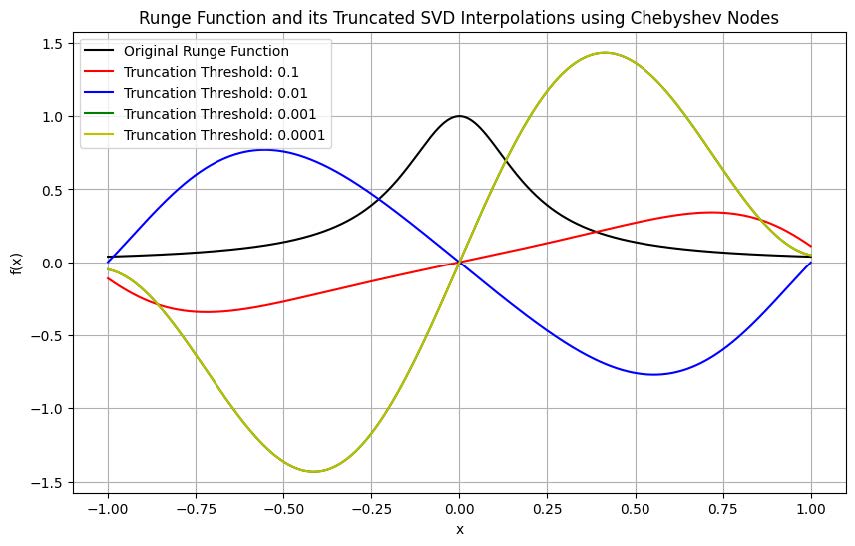}
    \caption{Runge Function and its Truncated SVD Interpolations using Chebyshev Nodes}
    \label{fig:svd_chebyshev_nodes}
\end{figure}

Visualizations depict the real Runge function alongside the interpolation outcomes at various truncation thresholds. Figures \ref{fig:svd_truncation}, \ref{fig:svd_uniform_grid}, and \ref{fig:svd_chebyshev_nodes} illustrate the comparative smoothness and precision of each method. However, suboptimal results at the endpoints suggest revisiting the interpolation with Chebyshev nodes, potentially yielding improvements due to their denser distribution, especially at interval extremes.

\section{Conclusion and Outlook}
This paper rigorously examines various interpolation methodologies aimed at mitigating the Runge phenomenon, a prevalent issue in processing certain oscillatory functions. Traditional numerical approaches, such as the employment of Chebyshev nodes and cubic splines, have demonstrated considerable effectiveness. Significantly, the incorporation of regularization techniques, which utilize statistical insights, has markedly improved interpolation outcomes over non-regularized methods. The development of mock-Chebyshev interpolation represents an enhancement of traditional Chebyshev interpolation, yielding greater practical utility.

Further advancements are noted in the TISI three-interval interpolation method, which effectively addresses the Runge phenomenon in oscillatory function interpolation. By integrating diverse interpolation techniques, this method achieves more refined and accurate results. Particularly, the enhanced three-interval interpolation method excels in handling the Runge function, notably at the endpoints, and employs Chebyshev nodes for superior interpolation in the central interval. Additionally, the constrained mock-Chebyshev least squares approach, combining mock-Chebyshev node selection with constrained least squares interpolation, offers another viable solution to the Runge phenomenon.

Looking ahead, the paper identifies potential areas for optimization and refinement in these methods. For instance, the TISI method could benefit from exploring alternative interval division strategies or employing varied interpolation techniques within each interval. In the case of SVD-based interpolation, further research is needed to ascertain the optimal truncation point and enhance interpolation across different node types. Future studies might incorporate additional constraints or integrate other mathematical tools, like functional analysis and optimization techniques, to further refine interpolation outcomes. Moreover, the exploration of novel interpolation strategies to address the Runge phenomenon and related challenges remains an open field. Ultimately, empirical validation with practical data sets is crucial to ascertain the effectiveness and reliability of these methods in real-world applications.

\end{document}